\documentclass[runningheads]{llncs}
\usepackage[T1]{fontenc}
\usepackage{graphicx}
\usepackage{amsmath} 
\usepackage{amsfonts} 
\usepackage{float} 
\usepackage{comment} 
\begin{document}
\title{Solving the n-Queens Problem in Higher Dimensions}
%
%
\author{Tim Kunt\inst{1}\orcidID{0009-0006-5732-3208}}
\authorrunning{T. Kunt}
%
\institute{Zuse Institute, Takustraße 7, 14195 Berlin, Germany \\
\email{kunt@zib.de}, 
\url{kunt.org}}
%
\maketitle              
\begin{abstract}
How many mutually non-attacking queens can be placed on a d-dimensional chessboard of size n? The n-queens problem in higher dimensions is a generalization of the well-known n-queens problem. 
We present an integer programming formulation of the n-queens problem in higher dimensions and several strengthenings through additional valid inequalities. Compared to recent benchmarks, we achieve a speedup in computational time between 15-70x over all instances of the integer programs. Our computational results prove optimality of certificates for several large instances. Breaking additional, previously unsolved instances with the proposed methods is likely possible. On the primal side, we further discuss heuristic approaches to constructing solutions that turn out to be optimal when compared to the IP. 

\end{abstract}

\section{Problem Definition}
The classical $n$-queens problem is concerned with the question of how to place $n$ queens on an $n \times n$ chessboard, such that no pair of queens attack each other, i.e. that they do not share a row, column or diagonal. For the $n$-queens problem in higher dimensions, we consider a board of size $n$ and dimension $d$, where integer $n$ is the number of squares in each dimension and integer $d$ is the number of dimensions the board extends to. 
\begin{definition}[Queen] A queen \textbf{q} on a $(n,d)$-board is denoted as a tuple of integers $(a_1, a_2, ... , a_d)$ with $1 \leq a_i \leq n$ for $i \in \{1,2,...,d\}$.\\
For any pair of queens $q_1 = (a_1, a_2, ... , a_d)$ and $q_2 = (b_1, b_2, ... , b_d)$ with $q_1 \neq q_2$, we say that \textbf{$q_1$ attacks $q_2$} if there exists some $m\in \mathbb{Z}$ and nonzero  $\epsilon = (\epsilon_1, \epsilon_2, ... , \epsilon_d)$ with $\epsilon_i \in \{-1,0,1\}$ such that
\begin{align}
    a_i = \epsilon_i \cdot m + b_i \; \; \forall i \in \{1,2,...,d\}
\end{align}
\end{definition}
\textbf{Problem 1 ($(n,d)$-queens).} \vspace{0.2cm} \\
\begin{tabular}{rll}
    \textsc{Problem:} & \textit{A $(n,d)$-board,} $n,d \in \mathbb{N}^{+}$\vspace{0.2cm} \\
    \textsc{Solution:} & \textit{A set $Q$ of queens on the $(n,d)$-board such that: \vspace{0.0cm}}\\
    &- $|Q|=n^{d-1}$ \\
    &- \textit{For any distinct pair of queens} $q_1$, $q_2$ $\in Q$, $q_1$ \textit{does not attack} $q_2$ \vspace{0.2cm}
\end{tabular}\\
The above is a generalization of the definition \cite{gent2017complexity} of the standard $n$-queens problem to higher dimensions. While the $(n,2)$-queens problem has a solution for all $n \geq 3$ \cite{pauls1874maximalproblem}, a comparable result is not known for the $(n,d)$-queens problem for $d \geq 3$. For $n, d$, for which no solution exists, we are interested in finding sets of mutually non-attacking queens with maximal cardinality. This motivates the following definition of the \textit{partial} $(n,d)$-queens problem.\\
\\
\textbf{Problem 2 (partial $(n,d)$-queens).} \vspace{0.2cm} \\
\begin{tabular}{rll}
    \textsc{Problem:} & \textit{A $(n,d)$-board,} $k\leq n^{d-1}$, $n,d,k \in \mathbb{N}^{+}$\vspace{0.2cm} \\
    \textsc{Solution:} & \textit{A set $Q$ of queens on the $(n,d)$-board such that: \vspace{0.0cm}}\\
    &- $|Q|=k$ \\
    &-  \textit{For any distinct pair of queens} $q_1$, $q_2$ $\in Q$, $q_1$ \textit{does not attack} $q_2$ \vspace{0.2cm}
\end{tabular}\\
\\
For given $n, d$, the corresponding optimization problem asks to find the maximal $k$, for which a solution to the partial $(n,d)$-queens problem exists. We call such a solution a \textbf{maximal partial solution} \cite{bell2009} and denote it as $Q_{max}(n,d)$. 
We distinguish between the $(n,d)$-queens problem and the partial $(n,d)$-queens problem, as several results and methods differ or only apply for $k=n^{(d-1)}$.\\
\\
The generalization of the $n$-queens problem to the third dimension is first proposed by \cite{mccarty1978queen}. 
An overview on variants of the $n$-queens problem, including the problem in higher dimensions, is given by \cite{bell2009}. Recently \cite{langlois2022complexity} succeeded in proving maximality for several new instances. For a more comprehensive literature review and further theoretical results on both lower and upper bounds
, we refer to the dissertation of the main author \cite{kunt2024n}.

\section{Integer Programming Formulation}
\subsection{Base Model}
This integer programming formulation considers the $\frac{(3^d-1)}{2}$ hyperplanes of attack (c.f. \cite{Nudelman1995}) of a queen on the $(n,d)$-board. It restricts the sum over all squares in each of those lines to be less or equal to $1$, as two or more queens in one such line would attack each other. 
Denote $S$ the set of all squares on the $(n,d)$-board and $H$ the set of all hyperplanes of attack. Let $L_h$ be the superset of all sets of squares that lie in one of the hyperplanes $h \in H$. Then the partial $(n,d)$-queens problem can be formulated as follows
\begin{alignat*}{3}
  & \text{max } &       & \sum_{s \in S}^{n} x_s \\[2ex]
  & \text{s.t } & \quad &\; \sum_{s \in L }^{n} x_s  \leq 1,  && \quad \forall L \in L_{h}, \quad  \forall h \in H 
\end{alignat*}

\subsection{Strengthening the IP}
We present a set of valid inequalities for the partial $(n,d)$-queens problem, which provide upper bounds and strengthen the IP by restricting the LP relaxation. 
\subsubsection{Subsolution Inequalities}
For any subset of the $(n,d)$-board, we may apply our knowledge of maximal partial solutions for $(m,d)$-boards, $m < n$. Let $S^{n,d}_m$ denote the set of all subsets of the $(n,d)$-board that correspond to $(m,d)$-boards, i.e., $m^d$ hypercubes. We obtain the following additional valid inequalities:
\begin{align}
     \sum_{s\in S} x_s \leq  |Q_{max}(m,d)| \;\;\;\;\; \forall S \in S^{n,d}_m
\end{align}
\subsubsection{Layer Inequalities}
Recall that $|Q_{max}(n,d)| \leq |Q_{max}(n,d-k)| \cdot n^k$. More specifically, for all layer subsets of the $(n,d)$-board, we may apply our knowledge of maximal partial solutions for $(n,d-1)$-boards. Let $L^{n,d}$ denote the set of all layers of the $(n,d)$-board, meaning all subsets that correspond to $(n,d-1)$-boards. We obtain the following additional valid inequalities:
\begin{align}
     \sum_{s\in L} x_s \leq  |Q_{max}(n,d-1)| \;\;\;\;\; \forall L \in L^{n,d}
\end{align}
These inequalities may be applied recursively for layers of layers. 
\subsubsection{Cube and Star Cliques}
\cite{fischetti2019finding} discuss a set of clique inequalities for the $(n,2)$-queens problem. These generalize nicely to higher dimensions, where they take the shape of hypercube and the cross-polytope in the respective dimension. 
For general $d$, we may describe the cube inequalities as follows.
\begin{figure}[H]
\setkeys{Gin}{width=0.49\linewidth}
\includegraphics{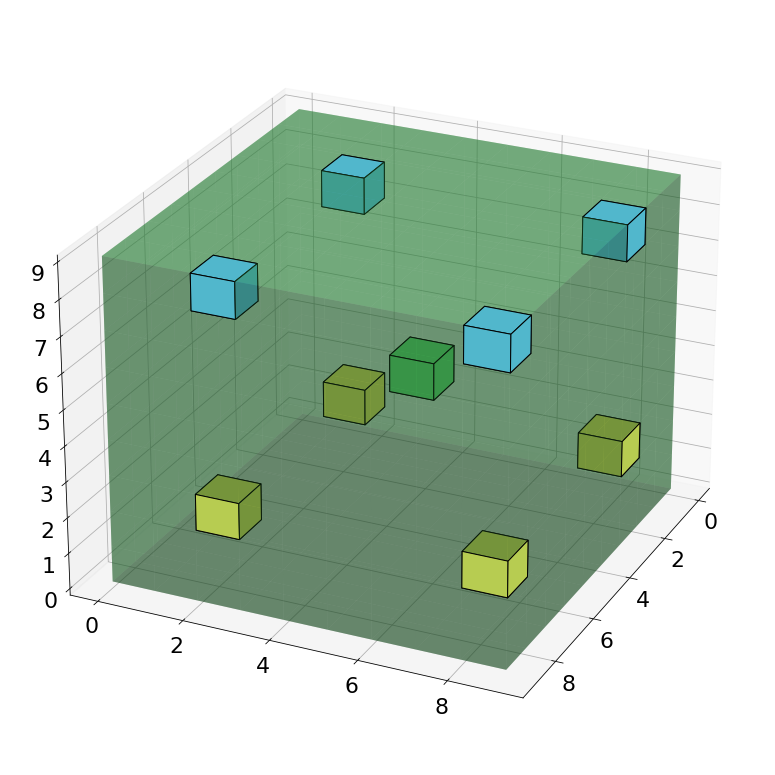}
\hfill
\includegraphics{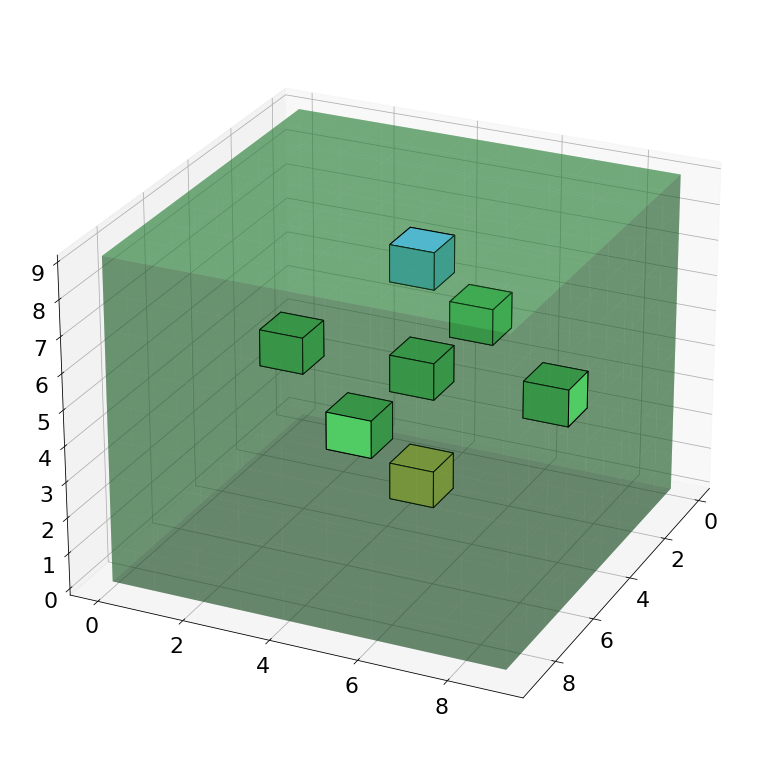}
\caption{Cube for $h=6$ and star clique for $h=3$, both on the $(9,3)$-board}
\label{cliques}
\end{figure}
\noindent For integer $h$ odd and $s_i+h \leq n$ for all $i \leq d$
\begin{align}
    & \bigg( \sum_{a_1 =0}^1 \sum_{a_2 =0}^1 \dots \sum_{a_d =0}^1 x_{(s_1+a_1 \cdot h,s_2+a_2 \cdot h,...,s_d+a_d \cdot h)} \bigg) \leq 1
\end{align}
And further for even $h$ 
\begin{align}
    & \Bigg( \bigg( \sum_{a_1 =0}^1 \sum_{a_2 =0}^1 \dots \sum_{a_d =0}^1 x_{(s_1+a_1 \cdot h,s_2+a_2 \cdot h,...,s_d+a_d \cdot h)} \bigg) \\
    & + x_{(s_1+\cdot \frac{h}{2},s_2+ \cdot \frac{h}{2},...,s_d+ \cdot \frac{h}{2})} \Bigg) \leq 1
\end{align}
For the star cliques we get for for $s_i+h \leq n$ and $s_i-h \geq n$ for all $i \leq d$
\begin{align}
    & x_{(s_1,s_2,...,s_d)} \nonumber \\
    & + x_{(s_1+h,s_2,...,s_d)} + x_{(s_1-h,s_2,...,s_d)} \nonumber \\
    & + x_{(s_1,s_2+h,...,s_d)} + x_{(s_1,s_2-h,...,s_d)} \nonumber \\
    & + ... \nonumber \\
    & + x_{(s_1,...,s_i+h,...,s_d)} + x_{(s_1,...,s_i-h,...,s_d)} \nonumber \\
    & + ... \nonumber \\
    & + x_{(s_1,s_2,...,s_d+h)} + x_{(s_1,s_2,...,s_d-h)} \leq 1
\end{align}
\section{Computational Results}
All instances were solved using Gurobi version 9.5.2 with 10 randomly generated seeds on a Intel Xeon Gold 6246 with 24 cores, 48 threads, and 384GB RAM. We compare our following different variants of the IP and the formulation by \cite{langlois2022complexity} regarding average runtime:
\begin{itemize}
    \item \textbf{Base} base model, only contains the necessary inequalities as described
    \item \textbf{LMR} the IP formulation by \cite{langlois2022complexity}, see \cite{Langlois2023}
    \item \textbf{Cube} base model with the additional cube clique inequalities
    \item \textbf{Star} base model with the additional star clique inequalities
    \item \textbf{C+S} base model with the additional cube- and star clique inequalities
    \item \textbf{Inf} proving infeasibility for $|Q_{max}|+1$, provided we know a certificate with $|Q_{max}|=k$, this allows us to prove that it is a maximal partial solution
    \item \textbf{WS} warmstarts, the solver is warmstarted with a certificate of a maximal partial solution
    \item \textbf{All} (only for $d\geq4$) base model with cube- and star cliques, layer- and subsolution inequalities
    \item further variants are combinations of the above
\end{itemize}
The fastest variant is marked bold for instances where the base model exceeds $10s$ runtime. Timeout is set at $10h$. Aim of the comparison is to understand which additional inequalities provide a speedup and subsequently which formulations may be able to break larger unsolved instances. Variance of runtimes differs among the instances, however the coefficient of variation is less than $1.0$ for all instances. Additionally, we ran larger instances that exceeded runtimes feasible for the above setup (and are therefore excluded to retain a fair comparison) on a PowerEdge C6520 with two Intel Xeon Gold 6338 with 64 cores, 128 threads, and 1TB RAM in total. This includes instances $d=3$, $n=9, 10$, which terminated in $31526s$ and $87914s$ respectively.
\begin{table}[H]
\noindent
\begin{tabular}{|rr|rrrrrrrrr|}
    \hline 
    $n$ & $|Q_{max}|$& Base   & LMR     & Cube     & Star    & C+S    & Inf    & Inf C+S & ws & ws C+S \\ \hline
    1 &  1           &$0.00$  &$0.00$   &$0.00$    &$0.00$   &$0.00$  &$0.00$  &$0.00$   & $0.00$ & $0.00$ \\ 
    2 &  1           &$0.00$  &$0.00$   &$0.00$    &$0.00$   &$0.00$  &$0.00$  &$0.00$   & $0.00$ & $0.00$ \\ 
    3 &  4           &$0.00$  &$0.00$   &$0.00$    &$0.00$   &$0.00$  &$0.00$  &$0.00$   & $0.00$ & $0.00$ \\ 
    4 &  7           &$1.03$  &$1.20$   &$0.81$    &$0.76$   &$0.81$  &$0.03$  &$0.42$   & $1.08$ & $0.81$ \\ 
    5 & 13           &$3.66$  &$3.25$   &$1.61$    &$3.48$   &$1.58$  &$2.48$  &$1.42$   & $4.30$ & $1.77$ \\ 
    6 & 21           &$135.71$&$101.12$ &$15.44$   &$163.24$ &$16.43$ &$349.06$&$\textbf{6.52}$   & $143.08$ & $18.14$ \\ 
    7 & 32           &$3351.14$&$6225.06$&$574.01$  &$3714.34$&$531.74$&$-$     &$\textbf{267.05}$ & $3839.39$ & $498.39$ \\ 
    8 & 48           &$-$     &$-$      &$11343.43$&$-$      &$11745.36$  &$-$ &$\textbf{4127.32}$&$-$ & $7324.73$ \\ 
    9 & 67           &$-$     &$-$      &$-$       &$-$      &$-$     &$-$     &$-$      & $-$  &$-$   \\ 
   10 & 91           &$-$     &$-$      &$-$       &$-$      &$-$     &$-$     &$-$      & $-$  &$-$   \\ 
   11 &121           &$894.79$&$2029.97$&$1164.61$ &$1246.48$     &$2001.73$&$\textbf{0.00}$  &$0.02$   &$1.10$&$6.70$\\ 
   12 &133           &$-$     &$-$      &$-$       &$-$      &$-$     &$-$     &$-$      & $-$  &$-$   \\ 
   13 &169           &$-$     &$-$      &$-$       &$-$      &$-$     &$\textbf{0.01}$  &$0.06$   &$2.74$&$30.66$\\ 
   14 &$\geq172$     &$-$     &$-$      &$-$       &$-$      &$-$     &$-$     &$-$      & $-$  &$-$   \\ \hline
\end{tabular}
\caption{\label{table-d3-runtimes}Average runtime comparison for $(n,3)$}
\end{table}
\vspace{-14mm}
\begin{table}[H]
\noindent
\begin{tabular}{|rr|rrrrrrr|}
    \hline 
    $n$ & $|Q_{max}|$& Base   & LMR     & All   & Inf & Inf All  & ws    & ws All \\ \hline
    1 &  1           &$0.00$  &$0.00$   &$0.00$  &$0.00$  &$0.00$   &$0.00$  &$0.00$   \\ 
    2 &  1           &$0.00$  &$0.00$   &$0.00$  &$0.00$  &$0.00$   &$0.00$  &$0.00$   \\ 
    3 &  6           &$0.02$  &$0.01$   &$0.02$  &$0.02$  &$0.02$   &$0.02$  &$0.02$   \\ 
    4 & 16           &$0.80$  &$0.09$   &$0.21$  &$0.67$  &$0.01$   &$0.57$  &$0.12$   \\ 
    5 & 38           &$-$     &$-$      &$-$     & $-$    & $-$      &$-$      &$-$    \\ 
    6 & 80           &$569.57$&$760.08$ &$725.91$&$879.43$&$\textbf{10.68}$&$64.98$  &$19.77$   \\ 
    7 &145           &\hphantom{yyyyyyy}$-$     &\hphantom{yyyyyyy}$-$&\hphantom{yyyyyyy}$-$      &\hphantom{yyyyyyy}$-$      &\hphantom{yyyyyyy}$-$        &\hphantom{yyyyyyy}$-$      &\hphantom{yyyyyyy}$-$     \\ \hline
\end{tabular}
\caption{\label{table-d4-runtimes}Average runtime comparison for $(n,4)$}
\end{table}
\vspace{-14mm}
\begin{table}[H]
\noindent
\begin{tabular}{|rr|rrrrr|}
    \hline 
    $n$ & $|Q_{max}|$& Base   & LMR     & All   & Inf & Inf All  \\ \hline
    1 &  1           &$0.00$  &$0.00$   &$0.00$  &$0.00$  &$0.00$      \\ 
    2 &  1           &$0.00$  &$0.00$   &$0.00$  &$0.00$  &$0.00$      \\ 
    3 & 11           &$0.78$  &$0.20$   &$0.53$  &$2.11$  &$0.31$     \\ 
    4 & 32           &$21.14$ &$3.38$   &$2.40$  &$44.47$  &$\textbf{0.07}$      \\ 
    5 &$\geq 89$           &\hphantom{yyyyyyy}$-$     &\hphantom{yyyyyyy}$-$&\hphantom{yyyyyyy}$-$      &\hphantom{yyyyyyy}$-$      &\hphantom{yyyyyyy}$-$  \\ \hline
\end{tabular}
\caption{\label{table-d4-runtimes}Average runtime comparison for $(n,5)$}
\end{table}
\noindent Constructing solutions with $|Q_{max}|$ queens for the discussed instances can be achieved with heuristic methods \cite{kunt2024n,Missouri}, their maximality however remains to be proven. Likewise the computationally expensive part of the IP is closing the dual. This is significantly improved by both introducing the discussed inequalities and by proving infeasibility for solutions of greater size than the known certificates instead. For $d=3$, $n=11,13$, the infeasibility of a larger solution $|Q_{max}|+1$ is trivially solved as it already violates the LP relaxation\footnote{Note that if $\text{gcd}(n,(2^d-1)!) = 1$ then $n^{d-1}$ non-attacking queens can be placed on the $(n,d)$-board, i.e. the $(n,d)$-queens problem has a solution \cite{van1981latin,klarner1979queen}}. \\
We emphasize that this problem  differs from the classic $(n,2)$-queens problem in complexity and, more importantly in structure, as the latter comes down to a feasibility problem and does not inherit the difficulty of closing the dual. Here, we showed that certain clique inequalities, in particular, the (hyper-)cube inequalities as discussed by \cite{fischetti2019finding} for the $(n,2)$-queens problem improve on the dual and achieve a speedup of $15.5\times -71.2 \times $ less computational time compared to \cite{langlois2022complexity}, who recently succeeded in breaking new, never-before-solved instances.\\
Following this conclusion, a recipe for solving further instances seems to be a combination of  (a) further improving primal heuristics, either through construction methods or different approaches such as QUBO, (b) study of clique inequalities and clique separation, and (c) making use of the unique structure of the problem during the solving process, for example through its corresponding symmetry group. Going forward, we also plan to compare the IP with constraint programming and SAT solvers, specifically for those instances where we suspect we already know the maximal partial solution.
%
%
\subsubsection{\discintname}
The author has no competing interests to declare that are
relevant to the content of this article.
\bibliographystyle{splncs04}
\bibliography{main}
%
\end{document}